\documentclass{autart}
\usepackage{natbib}

\usepackage{amsmath,amssymb,theorem}
\usepackage{graphicx}
\usepackage{algorithm,algorithmic}
\usepackage{psfrag}
	
\usepackage{color,xspace}
\usepackage{subfigure}

\newtheorem{theorem}{Theorem}[section]

\newtheorem{remark}[theorem]{Remark}
	
\newtheorem{corollary}[theorem]{Corollary}
\newtheorem{proposition}[theorem]{Proposition}

\newcommand{\real}{{\mathbb{R}}}
\newcommand{\realpositive}{\mathbb{R}_{>0}}
\newcommand{\realnonnegative}{\mathbb{R}_{\ge 0}}
\newcommand{\integernonnegative}{\mathbb{Z}_{\ge 0}}
\newcommand{\integerpositive}{\mathbb{Z}_{> 0}}

\newcommand{\GG}{{\mathcal{G}}}

\newcommand{\douti}{d_i^{\operatorname{out}}}
\newcommand{\Nouti}{\mathcal{N}_i^{\operatorname{out}}}
\newcommand{\Noutmax}{\mathcal{N}_{\max}^{\operatorname{out}}}
\newcommand{\Nini}{\mathcal{N}_i^{\operatorname{in}}}
\newcommand{\Dout}{D^{\operatorname{out}}}
\newcommand{\Din}{D^{\operatorname{in}}}

\renewcommand{\AA}{{\mathcal{A}}}
\newcommand{\BB}{{\mathcal{B}}}

\newcommand{\ones}[1]{\mathbf{1}_{#1}}
\newcommand{\zeros}[1]{\mathbf{0}_{#1}}

\newcommand{\diag}[1]{\operatorname{diag}\left( #1\right)}

\newcommand{\eps}{\varepsilon}
\renewcommand{\epsilon}{\varepsilon}

\renewcommand{\time}[1]{t_{\text{#1}}}
\newcommand{\Norm}[1]{\|#1\|}
\newcommand{\Sym}[1]{{#1}_{s}}

\newcommand{\until}[1]{\{1,\dots, #1\}}

\newcommand{\setdef}[2]{\{#1 \; | \; #2\}}

\renewcommand{\hat}{\widehat}

\newcommand{\commgraph}{\GG_\text{comm}}

\renewcommand{\commgraph}{\GG}

\newcommand{\oprocendsymbol}{\hbox{$\bullet$}}
\newcommand{\oprocend}{\relax\ifmmode\else\unskip\hfill\fi\oprocendsymbol}

\newcommand{\longthmtitle}[1]{\mbox{}\textup{\textbf{(#1)}}}

\renewcommand{\theenumi}{(\roman{enumi})}

\newcommand{\algoeventconsensus}{\textsc{event-triggered communication and control law}\xspace} 
\newcommand{\algoperiodic}{\textsc{periodic event-triggered communication and control law}\xspace}

\begin{document}

\date{\today}

\runauthor{Cameron Nowzari and Jorge Cort\'es}

\begin{frontmatter}

  \title{Distributed event-triggered coordination for average
    consensus on weight-balanced digraphs}

  \thanks{A preliminary version was presented as~\citep{CN-JC:14-acc}
    at the 2014 American Control Conference.}

  \author[penn]{Cameron Nowzari} \ead{cnowzari@seas.upenn.edu}\qquad
  \author[pepe]{Jorge Cort\'es} \ead{cortes@ucsd.edu}
  
  \address[penn]{Department of Electrical and Systems Engineering, 
  University of Pennsylvania, Philadelphia, PA, 19104, USA}
  \address[pepe]{Department of Mechanical and Aerospace Engineering,
    University of California, San Diego, CA, 92093, USA}

  \begin{abstract}
    This paper proposes a novel distributed event-triggered
    algorithmic solution to the multi-agent average consensus problem
    for networks whose communication topology is described by
    weight-balanced, strongly connected digraphs.  The proposed
    event-triggered communication and control strategy does not rely
    on individual agents having continuous or periodic access to
    information about the state of their neighbors. In addition, it
    does not require the agents to have a priori knowledge of any
    global parameter to execute the algorithm.  We show that, under
    the proposed law, events cannot be triggered an infinite number of
    times in any finite period (i.e., no Zeno behavior), and that the
    resulting network executions provably converge to the average of
    the initial agents' states exponentially fast. We also provide
    weaker conditions on connectivity under which convergence is
    guaranteed when the communication topology is switching.  Finally,
    we also propose and analyze a periodic implementation of our
    algorithm where the relevant triggering functions do not need to
    be evaluated continuously.  Simulations illustrate our results and
    provide comparisons with other existing algorithms.
  \end{abstract}


\end{frontmatter}

\section{Introduction}

This paper studies the multi-agent average consensus problem, where a
group of agents seek to agree on the average of their initial
states. Due to its numerous applications in networked systems, many
algorithmic solutions exist to this problem; however, a majority of
them rely on agents having continuous or periodic availability of
information from other agents.  Unfortunately, this assumption leads
to inefficient implementations in terms of energy consumption,
communication bandwidth, congestion, and processor usage.  Motivated
by these observations, our main goal here is the design of a provably
correct distributed event-triggered strategy that prescribes when
communication and control updates should occur so that the resulting
asynchronous network executions still achieve average consensus.

\emph{Literature review:} Triggered control seeks to understand the
trade-offs between computation, communication, sensing, and actuator
effort in achieving a desired task with a guaranteed level of
performance.  Early works~\citep{KJA-BMB:02} consider tuning
controller executions to the state evolution of a given system, but
the ideas have since then been extended to consider other
tasks, 
see~\citep{WPMHH-KHJ-PT:12} and references therein for a recent
overview.  Among the many references in the context of multi-agent
systems,~\citep{MMJ-PT:10} specifies the responsibility of each agent
in updating the control signals,~\citep{XW-MDL:11} considers network
scenarios with disturbances, communication delays, and packet drops,
and~\citep{CS-DV-JL:13} studies decentralized event-based control that
incorporates estimators of the interconnection signals among agents.
Several works have explored the application of event-triggered ideas
to the acquisition of information by the agents.
To this end,~\cite{GX-HL-LW-YJ:09,WPMHH-MCFD:13,XM-TC:13} combine
event-triggered controller updates with sampled data that allows for
the periodic evaluation of the triggers.~\cite{MZ-CGC:10} drop the
need for periodic access to information by considering event-based
broadcasts, where agents decide with local information only when to
obtain further information about neighbors.  Self-triggered
control~\citep{AA-PT:10,XW-MDL:09} relaxes the need for local
information by deciding when a future sample of the state should be
taken based on the available information from the last sampled
state. Team-triggered coordination~\citep{CN-JC:16-tac} combines the
strengths of event- and self-triggered control into a unified approach
for networked systems.

The literature on multi-agent average consensus is vast, see
e.g.,~\citep{ROS-JAF-RMM:07,WR-RWB:08,MM-ME:10} and references
therein.  \cite{ROS-RMM:03c} introduce a continuous-time algorithm
that achieves asymptotic convergence to average consensus for both
undirected and weight-balanced directed graphs. \cite{DVD-EF-KHJ:12}
build on this algorithm to propose a Lyapunov-based event-triggered
strategy that dictates when agents should update their control signals
but its implementation relies on each agent having perfect information
about their neighbors at all times.  The work~\citep{GSS-DVD-KHJ:13}
uses event-triggered broadcasting with time-dependent triggering
functions to provide an algorithm where each agent only requires exact
information about itself, rather than its neighbors. However, its
implementation requires knowledge of the algebraic connectivity of the
network. In addition, the strictly time-dependent nature of the
thresholds makes the network executions decoupled from the actual
state of the agents.  Closer to our treatment here,
\cite{EG-YC-HY-PA-DC:13} propose an event-triggered broadcasting law
with state-dependent triggering functions where agents do not rely on
the availability of continuous information about their neighbors
(under the assumption that all agents have initial access to a common
parameter). This algorithm works for networks with undirected
communication topologies and guarantees that all inter-event times are
strictly positive, but does not discard the possibility of an infinite
number of events happening in a finite time period.  We consider here
a more general class of communication topologies described by
weight-balanced, directed graphs. The
works~\citep{BG-JC:09c,AR-TC-CNH:14} present provably correct
distributed strategies that, given a directed communication topology,
allow a network of agents to find such weight edge assignments.

\emph{Statement of contributions:} Our main contribution is the design
and analysis of novel event-triggered broadcasting and controller
update strategies to solve the multi-agent average consensus problem
over weight-balanced digraphs.
With respect to the conference version of this
work~\citep{CN-JC:14-acc}, the present manuscript introduces new
trigger designs, extends the treatment from undirected graphs to
weight-balanced digraphs, and provides a comprehensive technical
treatment.  Our proposed law does not require individual agents to
have continuous access to information about the state of their
neighbors and is fully distributed in the sense that it does not
require any a priori knowledge by agents of global network parameters
to execute the algorithm.  Our Lyapunov-based design builds on the
evolution of the network disagreement to synthesize triggers that
agents can evaluate using locally available information to make
decisions about when to broadcast their current state to neighbors.
In our design, we carefully take into account the discontinuities in
the information available to the agents caused by broadcasts received
from neighbors and their effect on the feasibility of the resulting
implementation.  Our analysis shows that the resulting asynchronous
network executions are free from Zeno behavior, i.e., only a finite
number of events are triggered in any finite time period, and
exponentially converge to agreement on the average of all agents'
initial states over weight-balanced, strongly connected digraphs.  We
also provide a lower bound on the exponential convergence rate and
characterize the asymptotic convergence of the network under switching
topologies that remain weight-balanced and are jointly strongly
connected. Lastly, we propose a periodic implementation of our
event-triggered design that has agents check the triggers periodically
and characterize the sampling period that guarantees
correctness. Various simulations illustrate our results.


\section{Preliminaries}\label{se:preliminaries}

This section introduces some notational conventions and notions on
graph theory.  Let $\real$, $\realpositive$, $\realnonnegative$, and
$\integerpositive$ denote the set of real, positive real, nonnegative
real, and positive integer numbers, respectively. We denote by
$\ones{N}$ and $\zeros{N} \in \real^N$ the column vectors with entries
all equal to one and zero, respectively. We let $\Norm{\cdot}$ denote
the Euclidean norm on $\real^N$. We let $\diag{\real^N} = \setdef{x
  \in \real^N}{x_1=\dots=x_N}$.  For a finite set $S$, we let $|S|$
denote its cardinality.
Given $x, y \in \real$, Young's inequality~\citep{GHH-JEL-GP:52}
states that, for any $\epsilon \in \realpositive$,
\begin{align}\label{eq:Young}
  xy \leq \frac{x^2}{2\epsilon} + \frac{ \epsilon y^2}{2}.
\end{align}
A weighted directed graph (or weighted digraph) $\commgraph = (V,E,W)$
is comprised of a set of vertices $V = \until{N}$, directed edges $E
\subset V \times V$ and weighted adjacency matrix $W \in
\realnonnegative^{N \times N}$. Given an edge $(i,j) \in E$, we refer
to $j$ as an out-neighbor of $i$ and $i$ as an in-neighbor of $j$.
The sets of out- and in-neighbors of a given node $i$ are $\Nouti$ and
$\Nini$, respectively. The weighted adjacency matrix $W \in \real^{N
  \times N}$ satisfies $w_{ij} > 0$ if $(i, j) \in E$ and $w_{ij} = 0$
otherwise.  A path from vertex $i$ to $j$ is an ordered sequence of
vertices such that each intermediate pair of vertices is an edge.  A
digraph $\commgraph$ is strongly connected if there exists a path from
all $i \in V$ to all $j \in V$.  The out- and in-degree matrices
$\Dout$ and $\Din$ are diagonal matrices where
\begin{align*}
  \douti = \sum_{j \in \Nouti} w_{ij} , \quad d_i^{\operatorname{in}} =
  \sum_{j \in \Nini} w_{ji} ,
\end{align*}
respectively. A digraph is weight-balanced if $\Dout = \Din$.  The
(weighted) Laplacian matrix is $L = \Dout - W$. Based on the structure
of $L$, at least one of its eigenvalues is zero and the rest of them
have nonnegative real parts.  If the digraph $\commgraph$ is strongly
connected, $0$ is simple with associated eigenvector $\ones{N}$.  The
digraph~$\commgraph$ is weight-balanced if and only if $\ones{N}^T L =
\zeros{N}$ if and only if $ \Sym{L}=\frac{1}{2}(L + L^T)$ is positive
semidefinite.  For a strongly connected and weight-balanced digraph,
zero is a simple eigenvalue of $\Sym{L} $.  In this case, we order its
eigenvalues as $\lambda_1 = 0<\lambda_2\leq \dots \leq \lambda_N$, and
note the inequality
\begin{align}\label{eq:LapBound}
  x^T L x \ge \lambda_2(\Sym{L}) \| x -\frac{1}{N} (\ones{N}^T x)
  \ones{N} \|^2 ,
\end{align}
for all $x \in \real^N$.  The following property will also be of use
later,
\begin{align}\label{eq:lap2-bound}
  \lambda_2(\Sym{L}) x^T L x \leq x^T \Sym{L}^2 x \leq
  \lambda_N(\Sym{L}) x^T L x .
\end{align}
This can be seen by noting that $\Sym{L}$ is diagonalizable and
rewriting $\Sym{L} = S^{-1} D S$, where $D$ is a diagonal matrix
containing the eigenvalues of~$\Sym{L}$.

\section{Problem statement}\label{se:statement}

We consider the multi-agent average consensus problem for a network of
$N$ agents. We let $\commgraph$ denote the weight-balanced, strongly
connected digraph describing the communication topology of the
network.  Without loss of generality, we use the convention that an
agent $i$ is able to receive information from neighbors in $\Nouti$
and send information to neighbors in $ \Nini$.  We denote by $x_i \in
\real$ the state of agent $i \in \until{N}$. We consider
single-integrator dynamics
\begin{align}\label{eq:dynamics}
  \dot{x}_i(t) = u_i(t),
\end{align}
for all $i \in \until{N}$.  It is well known~\citep{ROS-RMM:03c} that
the distributed continuous control law
\begin{align}\label{eq:continuouslaw}
  u_i(t) = - \sum_{j \in \Nouti} w_{ij} (x_i(t) - x_j(t)) ,
\end{align}
drives each agent of the system to asymptotically converge to the
average of the agents' initial conditions. In compact form, this can
be expressed by
\begin{align*}
  \dot x(t) = -Lx(t) ,
\end{align*}
where $x(t) = (x_1(t), \dots, x_N(t))$ is the column vector of all
agent states and $L$ is the Laplacian of $\commgraph$.  However, in
order to be implemented, this control law requires each agent to
continuously access state information about its neighbors and
continuously update its control law. Here, we are interested in
controller implementations that relax both of these requirements by
having agents decide in an opportunistic fashion when to perform these
actions.

Under this framework, neighbors of a given agent only receive state
information from it when this agent decides to broadcast its state to
them. Equipped with this information, the neighbors update their
respective control laws.  We denote by $\hat{x}_i(t)$ the last
broadcast state of agent $i \in \until{N}$ at any given time~$t \in
\realnonnegative$.  We assume that each agent has continuous access to
its own state. We then utilize an event-triggered implementation of
the controller~\eqref{eq:continuouslaw} given by
\begin{align}\label{eq:control}
  u_i(t) = - \sum_{j \in \Nouti} w_{ij} (\hat{x}_i(t) - \hat{x}_j(t)).
\end{align}
Letting $u(t) = (u_1(t), \dots, u_N(t)) \in \real^N$ and 
$\hat{x} = ( \hat{x}_1, \dots, \hat{x}_N) \in \real^N$, we
write~\eqref{eq:control}~as
\begin{align*}
  u(t) = -L \hat{x} .
\end{align*}
Note that although agent $i$ has access to its own state $x_i(t)$, the
controller~\eqref{eq:control} uses the last broadcast state
$\hat{x}_i(t)$. This is to ensure that the average of the agents'
initial states is preserved throughout the evolution of the
system. More specifically, using this controller, one has
\begin{align}\label{eq:average-preservation}
  \frac{d}{dt} (\ones{N}^T x(t)) = \ones{N}^T \dot x(t) = -\ones{N}^T
  L \hat{x}(t) = 0 ,
\end{align}
where we
have used the fact that~$\commgraph$ is weight-balanced.

Our aim is to identify triggers that prescribe in an opportunistic
fashion when agents should broadcast their state to their neighbors so
that the network converges to the average of the initial agents'
states. Given that the average is conserved by~\eqref{eq:control}, all
the triggers should enforce is that the agents' states ultimately
agree.

\section{Distributed trigger design}\label{se:design}

In this section we synthesize a distributed triggering strategy that
prescribes when agents should broadcast state information and update
their control signals. Our design builds on the analysis of the
evolution of the network disagreement characterized by the following 
candidate Lyapunov function,
\begin{align}\label{eq:lyap}
  V(x) = \frac{1}{2}(x-\bar{x})^T (x-\bar{x}) ,
\end{align}
where $\bar{x} = \frac{1}{N} (\ones{N}^T x) \ones{N}$ corresponds to
agreement at the average of the states of all agents. The next result
characterizes a local condition for all agents in the network such
that this candidate Lyapunov function is monotonically nonincreasing.

\begin{proposition}\longthmtitle{Evolution of network
    disagreement}\label{pr:event}
  For  $i \in \until{N}$, let $a_i \in \realpositive$ and denote
  by $e_i(t) = \hat{x}_i(t) - x_i(t)$ the error between agent $i$'s
  last broadcast state and its current state at time $t \in
  \realnonnegative$.~Then,
  \begin{align*}
    \dot{V}(t) &\leq - \frac{1}{2} \sum_{i=1}^N \sum_{j \in \Nouti}
    w_{ij} \left[ (1-a_i)(\hat{x}_i-\hat{x}_j)^2 - \frac{e_i^2}{a_i}
    \right] .
  \end{align*}
\end{proposition}
\begin{pf}
  Note that, since the average is preserved,
  cf.~\eqref{eq:average-preservation}, under the control
  law~\eqref{eq:control}, $\bar{x} = \frac{1}{N} (\ones{N}^T x(0))
  \ones{N}$. The function $t \mapsto V(x(t))$ is continuous and
  piecewise continuously differentiable, with points of discontinuity
  of $\dot V$ corresponding to instants of time where an agent
  broadcasts its state. Whenever defined, this derivative takes the
  form
  \begin{align*}
    \dot{V} &= x^T \dot{x} - \bar{x}^T \dot{x} = -x^T L \hat{x} -
    \bar{x}^T L \hat{x} = -x^T L \hat{x},
  \end{align*}
  where we have used that the graph is weight-balanced in the last
  equality.  Let $e = (e_1, \dots, e_N) \in \real^N$ be the vector of
  errors of all agents. We can then rewrite $\dot{V}$ as
  \begin{align*}
    \dot{V} &= - \hat{x}^T L \hat{x} + e^T L \hat{x}.
  \end{align*}
  Expanding this out yields
  \begin{align*}
    \dot{V} = - \sum_{i=1}^N \sum_{j \in \Nouti} \left[ \frac{1}{2}
      w_{ij} (\hat{x}_i - \hat{x}_j)^2 - e_i w_{ij} (\hat{x}_i -
      \hat{x}_j) \right].
  \end{align*}
  Using Young's inequality~\eqref{eq:Young} for each product $e_i
  (\hat{x}_i - \hat{x}_j)$ with $\eps = a_i > 0$ yields,
  \begin{align*}
    \dot{V} &\leq - \sum_{i=1}^N \sum_{j \in \Nouti} w_{ij} \left[
      \frac{1}{2}(\hat{x}_i - \hat{x}_j)^2 - \frac{e_i^2}{2a_i} -
      \frac{a_i(\hat{x}_i - \hat{x}_j)^2}{2} \right]
    \\
    &= - \frac{1}{2} \sum_{i=1}^N \sum_{j \in \Nouti} w_{ij} \left[
      (1-a_i)(\hat{x}_i-\hat{x}_j)^2 - \frac{e_i^2}{a_i} \right] ,
  \end{align*}
  which concludes the proof.  \qed
\end{pf}

From Proposition~\ref{pr:event}, a sufficient condition to ensure that
the proposed candidate Lyapunov function $V$ is monotonically
decreasing is to maintain 
\begin{align*}
  \sum_{j \in \Nouti} w_{ij} \left[ (1-a_i)(\hat{x}_i-\hat{x}_j)^2 -
    \frac{e_i^2}{a_i} \right] \geq 0,
\end{align*}
for all $i \in \until{N}$ at all times.  This is accomplished by
ensuring
\begin{align*}
  e_i^2 \leq \frac{a_i(1-a_i)}{\douti} \sum_{j \in \Nouti}
  w_{ij} (\hat{x}_i - \hat{x}_j)^2 ,
\end{align*}
for all $i \in \until{N}$.  The maximum of the function $a_i(1-a_i)$
in the domain $(0,\infty)$ is attained at $a_i = \frac{1}{2}$, so we
have each agent select this value to optimize the trigger design.  As
a consequence of the above discussion, we have the following result.

\begin{corollary}\label{co:trigger}
  For each $i \in \until{N}$, let $\sigma_i \in (0,1)$ and define
  \begin{align}\label{eq:ogtrigger}
    f_i(e_i) = e_i^2 - \sigma_i \frac{1}{4 \douti} \sum_{j \in
      \Nouti} w_{ij} (\hat{x}_i - \hat{x}_j)^2 .
  \end{align}
  If each agent~$i$ enforces the condition $ f_i(e_i(t)) \le 0$ at all
  times, then
  \begin{align*}
    \dot{V}(t) \leq -\sum_{i=1}^N \frac{1-\sigma_i}{4} \sum_{j \in
      \Nouti} w_{ij} (\hat{x}_i - \hat{x}_j)^2.
  \end{align*}
  (Note that the latter quantity is strictly negative for all $\hat{x}
  \notin \diag{\real^N}$ because the graph is strongly connected).
\end{corollary}

For each $i \in \until{N}$, we refer to the function $f_i$ defined in
Corollary~\ref{co:trigger} as the \emph{triggering function} and to
the condition $f_i(e_i)=0$ as the \emph{trigger}.  Note that the
design parameter~$\sigma_i$ affects how flexible the trigger is: as
the value of $\sigma_i$ is selected closer to $1$, the trigger is
enabled less frequently at the cost of agent $i$ contributing less to
the decrease of the Lyapunov function.

An important observation is that, since the triggering function $f_i$
depends on the last broadcast states~$\hat{x}$, a broadcast from a
neighbor of $i$ might cause a discontinuity in the evaluation of
$f(e_i)$, where just before the update was received, $f_i(e_i) <0$,
and immediately after, $f_i(e_i)>0$. Such event would make agent $i$
miss the trigger. Thus, rather than prescribing agent $i \in
\until{N}$ to broadcast its state when $f_i(e_i) = 0$, we instead
define an event~by either
\begin{align}\label{eq:finaltrigger}
  f_i(e_i) > 0,
\end{align}
or
\begin{align}\label{eq:finaltrigger2}
  f_i(e_i) = 0 \quad \text{and} \quad \phi_i \neq 0,
\end{align}
where for convenience, we use the shorthand notation
\begin{align*}
  \phi_i = \sum_{j \in \Nouti} w_{ij} (\hat{x}_i - \hat{x}_j)^2 \in
  \realnonnegative .
\end{align*}
We note the useful equality $\sum_{i=1}^N \phi_i = 2\, \hat{x}^T L
\hat{x}$. The reasoning behind these triggers is the following.  The
inequality~\eqref{eq:finaltrigger} makes sure that the discontinuities
of $\phi_i$ do not make the agent miss an event.  The
trigger~\eqref{eq:finaltrigger2} makes sure that the agent is not
required to continuously broadcast its state to neighbors when its
last broadcast state is in agreement with the states received from
them.

The triggers~\eqref{eq:finaltrigger} and~\eqref{eq:finaltrigger2} are a
generalization of the ones proposed in~\citep{EG-YC-HY-PA-DC:13}.
However, it is unknown whether they are sufficient to exclude the
possibility of Zeno behavior in the resulting executions.  To address
this issue, we prescribe the following additional trigger.  Let
$\time{last}^i$ be the last time at which agent $i$ broadcast its
information to its neighbors~$\Nini$. If at some time $t \geq
\time{last}^i$, agent $i$ receives new information from a neighbor $j
\in \Nouti$, then $i$ immediately broadcasts its state if
\begin{align}\label{eq:designtrigger}
  t \in (\time{last}^i, \time{last}^i + \epsilon_i) .
\end{align}
Here, $\epsilon_i \in \realpositive$ is a design parameter selected
so~that
\begin{align}\label{eq:epsi}
  \epsilon_i < \sqrt{ \frac{\sigma_i}{4 \douti w^{\max}_i
      |\Nouti|} } ,
\end{align}
where $ w^{\max}_i = \max_{j \in \Nouti} w_{ij}$.  Our analysis in
Section~\ref{se:guarantees} will expand on the role of this bound and
the additional trigger in preventing the occurrence of Zeno behavior.

In conclusion, the
triggers~\eqref{eq:finaltrigger}-\eqref{eq:designtrigger} form the
basis of the \algoeventconsensus, which is formally presented in
Table~\ref{tab:algorithm}.

\begin{table}[htb]
  \centering
  \framebox[.9\linewidth]{\parbox{.85\linewidth}{%
      \parbox{\linewidth}{At all times $t$ agent $i \in \until{N}$
        performs:}
      \vspace*{-2.5ex}
      \begin{algorithmic}[1]
        \IF{
          $f_i(e_i(t)) > 0$ \textbf{or} ($f_i(e_i(t)) = 0$ and
          $\phi_i(t) \neq 0$)}  
        \STATE broadcast state information $x_i(t)$ and update control signal
        \ENDIF
        \IF{
          new information $x_j(t)$ is received from some neighbor(s)
          $j \in \Nouti$} 
        \IF{agent $i$ has broadcast its state at any time $t' \in (t - \epsilon_i, t)$} 
        \STATE broadcast state information $x_i(t)$ 
        \ENDIF
        \STATE update control signal
        \ENDIF
      \end{algorithmic}}}
  \caption{\algoeventconsensus.}\label{tab:algorithm}
\end{table}

Each time an event is triggered by an agent, say $i \in \until{N}$,
that agent broadcasts its current state to its out-neighbors and
updates its control signal, while its in-neighbors $j \in \Nini$
update their control signal.  This is in contrast to other
event-triggered designs, see e.g.,~\citep{LZ-CZ:10,DVD-EF-KHJ:12}, where
events only correspond to updates of control signals because exact
information is available to the agents at all times.


\section{Analysis of the \algoeventconsensus}\label{se:guarantees}

Here we analyze the properties of the control law~\eqref{eq:control}
in conjunction with the \algoeventconsensus of
Section~\ref{se:design}.
Our first result shows that the network executions are guaranteed not
to exhibit Zeno behavior. Its proof illustrates the role played by the
additional trigger~\eqref{eq:designtrigger} in facilitating the
analysis to establish this property.

\begin{proposition}\longthmtitle{No Zeno behavior}\label{prop:zeno}
  Given the system~\eqref{eq:dynamics} with control
  law~\eqref{eq:control} executing the \algoeventconsensus over a
  weight-balanced, strongly connected digraph, the agents will not be
  required to communicate an infinite number of times in any finite
  time period.
\end{proposition}
\begin{pf}
  We are interested in showing here that no agent will broadcast its
  state an infinite number of times in any finite time period.  Our
  first step consists of showing that, if an agent does not receive
  new information from neighbors, its inter-event times are lower
  bounded by a positive constant.
  Assume agent $i \in \until{N}$ has just broadcast its state at time
  $t_0$, and thus $e_i(t_0) = 0$. For $t \geq t_0$, while no new
  information is received, $\hat{x}_i(t)$ and $\hat{x}_j(t)$ remain
  constant.  Given that $\dot e_i = - \dot x_i$, the evolution of the
  error is simply
  \begin{align}\label{eq:error-evolution}
    e_i(t) = - (t - t_0) \hat{z}_i ,
  \end{align}
  where, for convenience, we use the shorthand notation $\hat{z}_i =
  \sum_{j \in \Nouti} w_{ij} (\hat{x}_j - \hat{x}_i)$.  Since we are
  considering the case when no neighbors of $i$ broadcast information,
  the trigger~\eqref{eq:designtrigger} is irrelevant.  We are then
  interested in finding the time~$t^*$ when~$f_i(e_i) = 0$ occurs,
  triggering a broadcast of agent $i$'s state.  If $\hat{z}_i = 0$, no
  broadcasts will ever happen ($t^* = \infty$) because $e_i(t) = 0$
  for all $t \geq t_0$.  Hence, consider the case when $\hat{z}_i \neq
  0$, which in turn implies $\phi_i \neq 0$.
  Using~\eqref{eq:error-evolution}, the
  trigger~\eqref{eq:finaltrigger2} prescribes a broadcast at the time
  $t^* \ge t_0$ satisfying
  \begin{align*}
    (t^* - t_0)^2
    \hat{z}_i^2 
    - \sigma_i \frac{1}{4 \douti} \sum_{j \in \Nouti} w_{ij}
    (\hat{x}_i - \hat{x}_j)^2 = 0,
  \end{align*}
  or, equivalently,
  \begin{align*}
    (t^* - t_0)^2 = \frac{ \sigma_i \sum_{j \in \Nouti} w_{ij}
      (\hat{x}_i - \hat{x}_j)^2 }{ 4 \douti \Big( \sum_{j \in
        \Nouti} w_{ij} (\hat{x}_i - \hat{x}_j) \Big)^2 }.
  \end{align*}
  Using the fact that $(\sum_{k=1}^p y_k)^2 \le p \sum_{k=1}^p y_k^2$
  for any $y_1,\dots,y_p \in \real$ and $p \in \integerpositive$
  (which readily follows from the Cauchy-Schwarz inequality), we
  obtain
  \begin{multline}\label{eq:great}
    \Big( \sum_{j \in \Nouti} w_{ij} (\hat{x}_i - \hat{x}_j) \Big)^2
    \leq |\Nouti| \sum_{j \in \Nouti} w_{ij}^2 (\hat{x}_i -
    \hat{x}_j)^2
    \\
    \leq |\Nouti| w^{\max}_i \sum_{j \in \Nouti} w_{ij} (\hat{x}_i -
    \hat{x}_j)^2 .
  \end{multline} 
  Therefore, we can lower bound the inter-event time by
  \begin{align*}
    \tau_i = t^* - t_0 \geq \sqrt{ \frac{\sigma_i}{4 \douti
        w^{\max}_i | \Nouti | }} > 0 ,
  \end{align*}
  (incidentally, this explains our choice in~\eqref{eq:epsi}).  Our
  second step builds on this fact to show that messages cannot be sent
  an infinite number of times between agents in a finite time period.
  Let time $t_0$ be the time at which agent $i$ has broadcast its
  information to neighbors and thus $e_i(t_0) = 0$. If no information
  is received by time $t_0 + \epsilon_i < t_0 + \tau_i$, there is no
  problem since $\epsilon_i > 0$, so we now consider the case that at
  least one neighbor of $i$ broadcasts its information at some time
  $t_1 \in (t_0, t_0 + \epsilon_i)$.  In this case it means that at
  least one neighbor $j \in \Nouti$ has broadcast new information,
  thus agent $i$ would also rebroadcast its information at time $t_1$
  due to trigger~\eqref{eq:designtrigger}.  Let $I$ denote the set of
  all agents who have broadcast information at time $t_1$ (we refer to
  these agents as synchronized). This means that, as long as no agent
  $k \notin I$ sends new information to any agent in $I$, the agents
  in $I$ will not broadcast new information for at least $\min_{j \in
    I} \tau_j$ seconds, which includes the original agent $i$. As
  before, if no new information is received by any agent in $I$ by
  time $t_1 + \min_{j \in I} \epsilon_j$ there is no problem, so we
  now consider the case that at least one agent $k$ sends new
  information to some agent $j \in I$ at time $t_2 \in (t_1, t_1 +
  \min_{j \in I} \epsilon_j)$. By trigger~\eqref{eq:designtrigger},
  this would require all agents in $I$ to also broadcast their state
  information at time $t_2$ and agent $k$ will now be added to the set
  $I$. Reasoning repeatedly in this way, the only way for infinite
  communications to occur in a finite time period is for an infinite
  number of agents to be added to the set $I$, which is not possible
  because there are only a finite number of agents~$N$. \qed
\end{pf} 



\begin{remark}[Conditions for Zeno] {\rm We note here that the
    introduction of the trigger~\eqref{eq:designtrigger} is sufficient
    to ensure Zeno behavior does not occur but it is an open problem
    to determine whether it is
    also necessary. 
    The design in~\cite[Corollary~2]{EG-YC-HY-PA-DC:13} specifies
    triggers of a nature similar
    to~\eqref{eq:finaltrigger}-\eqref{eq:finaltrigger2} for undirected
    graphs and guarantees that no agent undergoes an infinite number
    of updates at any given instant of time, but does not discard the
    possibility of an infinite number of updates in a finite time
    period, as Proposition~\ref{prop:zeno} does.  \oprocend }
\end{remark}

Next, we establish global exponential convergence.

\begin{theorem}\longthmtitle{Exponential convergence to average
    consensus}\label{th:exp-convergence}
  Given the system~\eqref{eq:dynamics} with control
  law~\eqref{eq:control} executing the \algoeventconsensus over a
  weight-balanced strongly connected digraph, all agents exponentially
  converge to the average of the initial states, i.e., $\lim_{t
    \rightarrow \infty} x(t) = \bar{x}$.
\end{theorem}
\begin{pf}
  By design, we know that the
  event-triggers~\eqref{eq:finaltrigger}-\eqref{eq:finaltrigger2}
  ensure that, cf. Corollary~\ref{co:trigger},
  \begin{align}\label{eq:Vdot}
    \dot{V} \leq \sum_{i=1}^N \frac{\sigma_i - 1}{4}
    \phi_i.
  \end{align}
  We show that convergence is exponential by establishing that the
  evolution of $V$ towards $0$ is exponential.  Define $\sigma_{\max}
  = \max_{i \in \until{N}} \sigma_i$ to further
  bounding~\eqref{eq:Vdot} by
  \begin{align*}
    \dot V \leq \frac{\sigma_{\max} - 1}{4}
    \sum_{i=1}^N \phi_i = \frac{\sigma_{\max} - 1}{2
      } \hat{x}^T L \hat{x} .
  \end{align*}
  Given this inequality, our next step is to relate the value of
  $V(x)$ with $\hat{x}^T L \hat{x}$.  Note that
  \begin{align*}
    V(x) & \le \frac{1}{2 \lambda_2(\Sym{L})} x^T L x = \frac{1}{2
      \lambda_2(\Sym{L})} (\hat{x} -e)^T L (\hat{x} -e)
    \\
    & = \frac{1}{2 \lambda_2(\Sym{L})} \big( \hat{x}^T L \hat{x} - 2
    \hat{x}^T \Sym{L} e + e^T L e \big),
  \end{align*}
  where we have used~\eqref{eq:LapBound} in the inequality. Now,
  \begin{align*}
    e^T L e \le \lambda_N(\Sym{L}) \| e \|^2 \le \lambda_N(\Sym{L})
    \frac{\sigma_{\max}}{2 d_{\min}^\text{out}} \hat{x}^T L \hat{x} ,
  \end{align*}
  where $d_{\min}^\text{out} = \min_{i \in \until{N}} \douti$
  and we have used $f_i(e_i)\le 0$ in the second inequality. On the
  other hand,
  \begin{align*}
    | \hat{x}^T \Sym{L} e | & \le \| \Sym{L} \hat{x}\| \, \|e\| \le
    \sqrt{\lambda_N(\Sym{L}) \hat{x}^T L \hat{x}} \,
    \sqrt{\frac{\sigma_{\max}}{2 d_{\min}^\text{out}} \hat{x}^T L
      \hat{x}}
    \\
    & = \sqrt{\lambda_N(\Sym{L}) \frac{\sigma_{\max} }{2
        d_{\min}^\text{out}}} \hat{x}^T L \hat{x},
  \end{align*}
  where we have used~\eqref{eq:lap2-bound} in the second
  inequality. Putting these bounds together, we obtain
  \begin{align*}
    V(x) & \le \AA \, \hat{x}^T L \hat{x},
  \end{align*}
  with $\AA =
  \tfrac{1}{2 \lambda_2(\Sym{L})} \big( 1 + \sqrt{\lambda_N(\Sym{L})
    \tfrac{\sigma_{\max} }{2 d_{\min}^\text{out}}} \big)^2$.
  Using this expression in the bound for the Lie derivative, we get
  \begin{align*}
    \dot V \leq \frac{\sigma_{\max} - 1}{2}
    \hat{x}^T L \hat{x} \le \frac{\sigma_{\max} - 1}{2
     \AA} V(x(t)).
  \end{align*}
  This, together with the fact that $t \mapsto V(x(t))$ is continuous
  and piecewise differentiable implies, using the Comparison Lemma,
  cf.~\citep{HKK:02}, that $V(x(t)) \le V(x(0)) \exp
  (\frac{\sigma_{\max} - 1}{2 \AA} t)$ and hence the exponential
  convergence of the network trajectories to the average state.  \qed
\end{pf}
%


The Lyapunov function used in the proof of
Theorem~\ref{th:exp-convergence} does not depend on the specific
network topology. Therefore, when the communication digraph is
time-varying, this function can be used as a common Lyapunov function
to establish asymptotic convergence to average consensus.  This
observation is key to establish the next result, whose proof we omit
for reasons of space.

\begin{proposition}\longthmtitle{Convergence under switching
    topologies}
  Let $\Xi_N$ be the set of weight-balanced digraphs over $N$
  vertices. Denote the communication digraph at time $t$ by $\GG(t)$.
  Consider the system~\eqref{eq:dynamics} with control
  law~\eqref{eq:control} executing the \algoeventconsensus over a
  switching digraph, where $t \mapsto \GG(t) \in \Xi_N$ is piecewise
  constant and such that there exists an infinite sequence of
  contiguous, nonempty, uniformly bounded time intervals over which
  the union of communication graphs is strongly connected. Then,
  assuming all agents are aware of who its neighbors are at each time
  and agents broadcast their state if their neighbors change, all
  agents asymptotically converge to the average of the initial states.
\end{proposition}

\section{Periodically checked event-triggered
  coordination}\label{se:periodic}

Here we propose an alternative strategy, termed \algoperiodic, where
agents only evaluate triggers~\eqref{eq:finaltrigger}
and~\eqref{eq:finaltrigger2} periodically, instead of continuously.
Specifically, given a sampling period $h \in \realpositive$, we let
$\{ t_\ell \}_{\ell \in \integernonnegative}$, where $t_{\ell+1} =
t_\ell + h$, denote the sequence of times at which agents evaluate the
decision of whether to broadcast their state to their neighbors. This
type of design is more in line with the constraints imposed by
real-time implementations, where individual components work at some
given frequency, rather than continuously. An inherent and convenient
feature of this strategy is the lack of Zeno behavior (since
inter-event times are naturally lower bounded by~$h$), making the need
for the additional trigger~\eqref{eq:designtrigger} superfluous. The
strategy is formally presented in Table~\ref{tab:algorithm2}.

\begin{table}[htb]
  \centering
  \framebox[.9\linewidth]{\parbox{.85\linewidth}{%
      \parbox{\linewidth}{At  times $t \in \{0, h, 2h, \dots \}$,
        agent $i \in \until{N}$ 
        performs:}
      \vspace*{-1.5ex}
      \begin{algorithmic}[1]
        \IF{
          $f_i(e_i(t)) > 0$ \textbf{or} ($f_i(e_i(t)) = 0$ and
          $\phi_i(t) \neq 0$)}  
        \STATE broadcast state information $x_i(t)$ and update control signal
        \ENDIF
        \IF{
          new information $x_j(t)$ is received from some neighbor(s)
          $j \in \Nouti$} 
        \STATE update control signal
        \ENDIF
      \end{algorithmic}}}
  \caption{\algoperiodic.}\label{tab:algorithm2}
\end{table}


Each time an agent~$i \in \until{N}$ broadcasts, this resets the error
to zero, $e_i = 0$. However, because triggers are not evaluated
continuously, we no longer have the guarantee $ f_i(e_i(t)) \leq 0$ at
all times $t$ but, instead, have
\begin{align}\label{eq:enforceperiodic}
  f_i(e_i(t_\ell)) \leq 0 ,
\end{align}
for $\ell \in \integernonnegative$.  The next result provides a
condition on~$h$ that guarantees the correctness of our design.

\begin{theorem}\longthmtitle{Exponential convergence under
    \algoperiodic} 
  Let $h \in \realpositive$ be such that
  \begin{align}\label{eq:bound}
    \sigma_{\max} + 4 h w_{\max} |\Noutmax | < 1,
  \end{align}
  where $ w_{\max} = \max_{i \in \until{N}} w_{i}^{\max}$ and
  $|\Noutmax | = \max_{i \in \until{N}} |\Nouti |$.
  Then, given the system~\eqref{eq:dynamics} with control
  law~\eqref{eq:control} executing the \algoperiodic over a
  weight-balanced strongly connected digraph, all agents exponentially
  converge to the average of the initial states.
\end{theorem}
\begin{pf}
  Since~\eqref{eq:enforceperiodic} is only guaranteed at the sampling
  times under the \algoperiodic, we analyze what happens to the
  Lyapunov function~$V$ in between them.  For $t \in [t_\ell,
  t_{\ell+1})$, note that
  \begin{align*}
    e(t) = e(t_\ell) + (t-t_\ell) L \hat{x}(t_\ell) .
  \end{align*}
  Substituting this expression into $\dot{V}(t) = -\hat{x}^T(t) L
  \hat{x}(t) + e^T(t) L \hat{x}(t)$, we obtain
  \begin{multline*}
    \dot{V}(t) = -\hat{x}^T(t_\ell) L \hat{x}(t_\ell) + e^T(t_\ell) L
    \hat{x}(t_\ell)
    \\
    + (t - t_\ell) \hat{x}^T(t_\ell) L^T L \hat{x}(t_\ell) ,
  \end{multline*}
  for all $t \in [t_\ell, t_{\ell+1})$.  For a simpler exposition,
  we drop all arguments referring to time $t_\ell$ in the sequel.
  Following the same line of reasoning as in
  Proposition~\ref{pr:event} yields
  \begin{align*}
    \dot{V}(t) \leq \sum_{i=1}^N \frac{\sigma_i - 1}{4} \phi_i +
    (t-t_\ell) \hat{x}^T L^T L \hat{x} .
  \end{align*}
  Using~\eqref{eq:great}, we bound
  \begin{align}\label{eq:boundxx}
    \hat{x}^T L^T L \hat{x} & = \sum_{i=1}^N \big(\sum_{j \in \Nouti}
    w_{ij} (\hat{x}_i-\hat{x}_j)\big)^2 \notag
    \\
    & \le \sum_{i=1}^N |\Nouti| w^{\max}_i \sum_{j \in \Nouti} w_{ij} (\hat{x}_i -
    \hat{x}_j)^2 \notag
    \\
    & = |\Noutmax | w_{\max} \sum_{i=1}^N \phi_i.
  \end{align}
  Hence, for $t \in [t_\ell, t_{\ell+1})$,
  \begin{align*}
    \dot{V}(t) & \leq \sum_{i=1}^N \Big( \frac{\sigma_i - 1}{4} + h
    w_{\max} |\Noutmax | \Big) \phi_i
    \\
    & \leq \Big( - \frac{1}{2} + \frac{\sigma_{\max}}{2} + 2 h
    w_{\max} |\Noutmax | \Big) \hat{x}^T L \hat{x} .
  \end{align*}
  Under~\eqref{eq:bound}, a reasoning similar to the proof of
  Theorem~\ref{th:exp-convergence} using~\eqref{eq:boundxx} leads to
  finding $\BB>0$ such~that
  \begin{align*}
    \dot{V}(t) & \leq \frac{1}{2 \BB }\Big( \sigma_{\max} + 4 h
    w_{\max} |\Noutmax | - 1 \Big) V(x(t)),
  \end{align*}
  which implies the result.  \qed
\end{pf}

\section{Simulations}\label{se:simulations}

This section illustrates the performance of the proposed algorithms in
simulation.  Figure~\ref{fig:sim1} shows a comparison of the
\algoeventconsensus with the algorithm proposed
in~\citep{EG-YC-HY-PA-DC:13} for undirected graphs over a network of
$5$ agents. Both algorithms operate under the
dynamics~\eqref{eq:dynamics} with control law~\eqref{eq:control}, and
differ in the way events are triggered. The algorithm
in~\citep{EG-YC-HY-PA-DC:13} requires all network agents to have
knowledge of an a priori chosen common parameter~$a \in
\realpositive$, which we set here to $a = 0.2$.
Figure~\ref{fig:sim1}(a) shows the evolution of the Lyapunov
function~$V$ and Figure~\ref{fig:sim1}(b) shows the number of events
triggered over time by each strategy.
%


{ \psfrag{oldwayonetwothree}[ll][ll]{\tiny \citep{EG-YC-HY-PA-DC:13}}%
  \psfrag{newwayonetwothree}[ll][ll]{\tiny Proposed}
  \psfrag{0}[cc][cc]{}
  \psfrag{1}[cc][cc]{\tiny $1$}%
  \psfrag{2}[cc][cc]{\tiny $2$} \psfrag{3}[cc][cc]{\tiny $3$}%
  \psfrag{4}[cc][cc]{\tiny $4$} \psfrag{5}[cc][cc]{\tiny $5$}%
  \psfrag{6}[cc][cc]{\tiny $6$} \psfrag{7}[cc][cc]{\tiny $7$}%
  \psfrag{8}[cc][cc]{\tiny $8$} \psfrag{10}[cc][cc]{\tiny $10$}
  \psfrag{12}[cc][cc]{\tiny $12$} \psfrag{14}[cc][cc]{\tiny $14$}
  \psfrag{16}[cc][cc]{\tiny $16$}
  \psfrag{3.5}[cc][cc]{\hspace*{-1.5mm}\tiny $3.5$}
  \psfrag{2.5}[cc][cc]{\hspace*{-1.5mm}\tiny $2.5$}
  \psfrag{1.5}[cc][cc]{\hspace*{-1.5mm}\tiny $1.5$}
  \psfrag{0.5}[cc][cc]{\hspace*{-1.5mm}\tiny $0.5$}
  \psfrag{50}[cc][cc]{\hspace*{-2mm}\tiny $50$}
  \psfrag{100}[cc][cc]{\hspace*{-2mm}\tiny $100$}
  \psfrag{150}[cc][cc]{\hspace*{-2mm}\tiny $150$}
  \begin{figure}[htb]
    \centering
    \subfigure[]{\includegraphics[width=.49\linewidth]{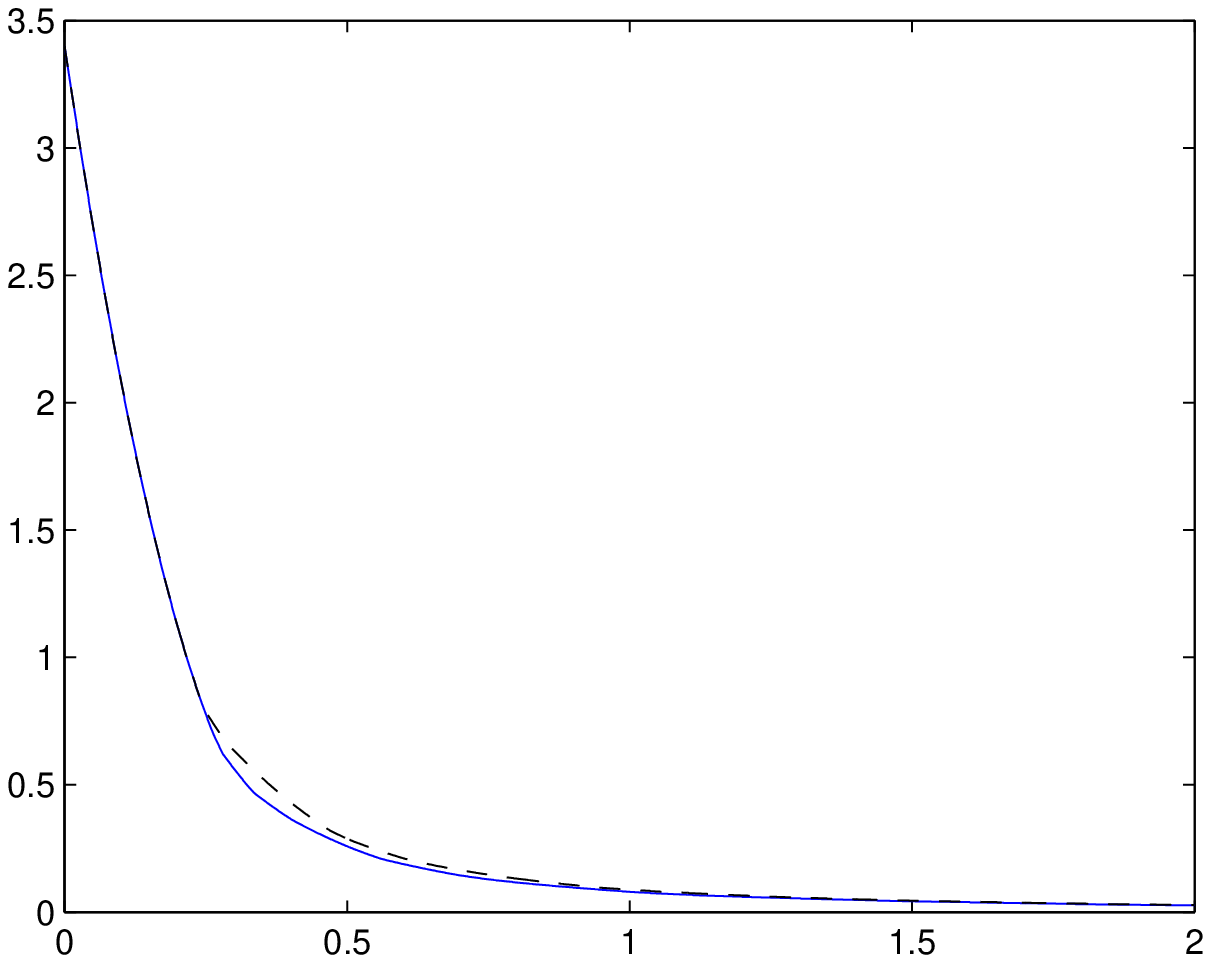}}
    \put(-104,5.5){\tiny $0$}
    \subfigure[]{\includegraphics[width=.49\linewidth]{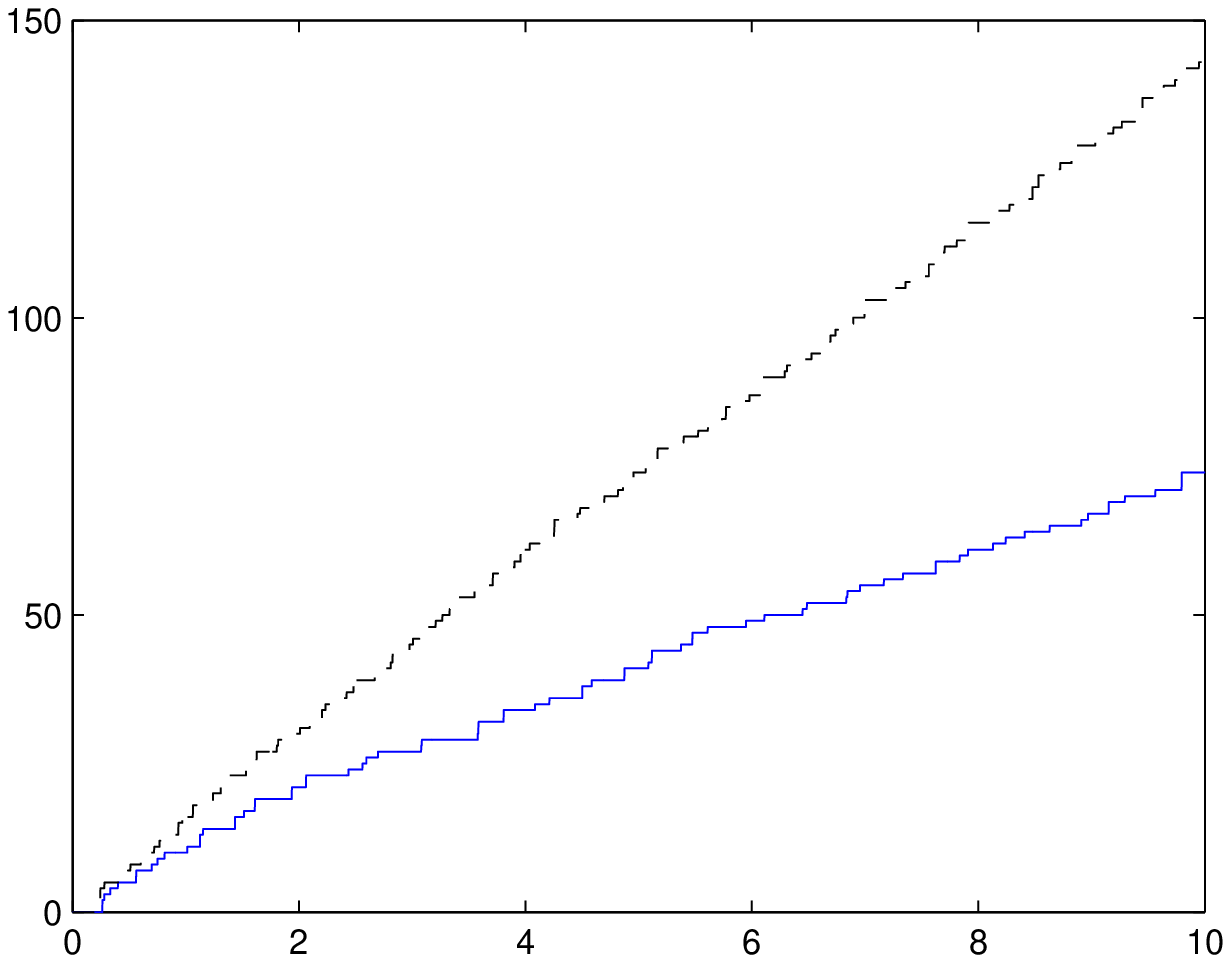}}
    \put(-65,0){\tiny Time}
    \put(-185,0){\tiny Time}
    \put(-235,68){\tiny $V$}
    \put(-120,65){\tiny $N_E$}
    \put(-104,5.5){\tiny $0$}
    \caption{Plots of (a) the evolution of the Lyapunov function $V$
      and (b) the total number $N_E$ of events of the
      \algoeventconsensus with $\sigma_i = 0.999$ for all~$i$ (solid blue) and the algorithm proposed
      in~\citep{EG-YC-HY-PA-DC:13} with $a = 0.2$ (dashed black). The
      network consists of $5$ agents with communication topology
      described by the undirected graph $(\until{5}, \{(1,2), (1,3),
      (2,4), (4,5)\})$.  The initial condition is $x(0) = [-1, 0, 2,
      2, 1]^T$.}\label{fig:sim1}
  \end{figure}
}

Figure~\ref{fig:sim2} shows an execution of \algoeventconsensus over a
network of $5$ agents whose communication topology is described by a
weight-balanced digraph. In this case, we do not compare it against
the algorithm in~\cite{EG-YC-HY-PA-DC:13} because the latter is only
designed to work for undirected graphs.

{ \psfrag{oldwayonetwothree}[ll][ll]{\tiny \citep{EG-YC-HY-PA-DC:13}}%
  \psfrag{newwayonetwothree}[ll][ll]{\tiny Proposed}
  \psfrag{0}[cc][cc]{}
  \psfrag{1}[cc][cc]{\tiny $1$}%
  \psfrag{2}[cc][cc]{\tiny $2$} \psfrag{3}[cc][cc]{\tiny $3$}%
  \psfrag{4}[cc][cc]{\tiny $4$} \psfrag{5}[cc][cc]{\tiny $5$}%
  \psfrag{6}[cc][cc]{\tiny $6$} \psfrag{7}[cc][cc]{\tiny $7$}%
  \psfrag{8}[cc][cc]{\tiny $8$} \psfrag{10}[cc][cc]{\hspace*{-2mm}\tiny $10$}
  \psfrag{12}[cc][cc]{\tiny $12$} \psfrag{14}[cc][cc]{\tiny $14$}
  \psfrag{16}[cc][cc]{\tiny $16$}
  \psfrag{3.5}[cc][cc]{}
  \psfrag{2.5}[cc][cc]{}
  \psfrag{1.5}[cc][cc]{}
  \psfrag{0.5}[cc][cc]{}
  \psfrag{5}[cc][cc]{\hspace*{-2mm}\tiny $5$}
  \psfrag{20}[cc][cc]{\hspace*{-2mm}\tiny $20$}
  \psfrag{40}[cc][cc]{\hspace*{-2mm}\tiny $40$}
  \psfrag{15}[cc][cc]{\hspace*{-2mm}\tiny $15$}
  \psfrag{35}[cc][cc]{\hspace*{-2mm}\tiny $35$}
  \psfrag{30}[cc][cc]{\hspace*{-2mm}\tiny $30$}
  \psfrag{25}[cc][cc]{\hspace*{-2mm}\tiny $25$}
  \psfrag{150}[cc][cc]{\hspace*{-2mm}\tiny $150$}
  \psfrag{200}[cc][cc]{\hspace*{-2mm}\tiny $200$}
  \psfrag{250}[cc][cc]{\hspace*{-2mm}\tiny $250$}
  \begin{figure}[htb]
    \centering
    \subfigure[]{\includegraphics[width=.49\linewidth]{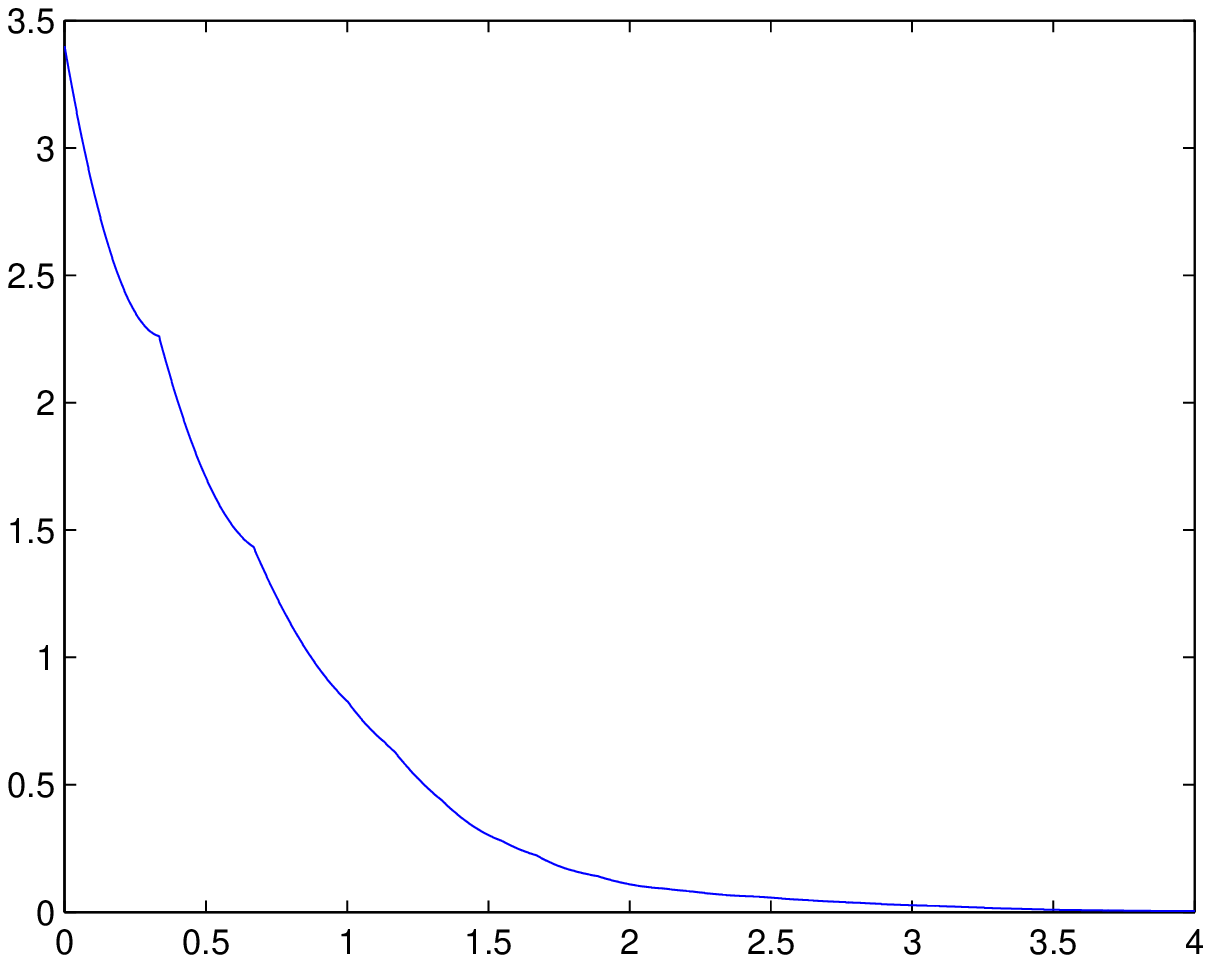}}
    \put(-104,5.5){\tiny $0$}
    \subfigure[]{\includegraphics[width=.49\linewidth]{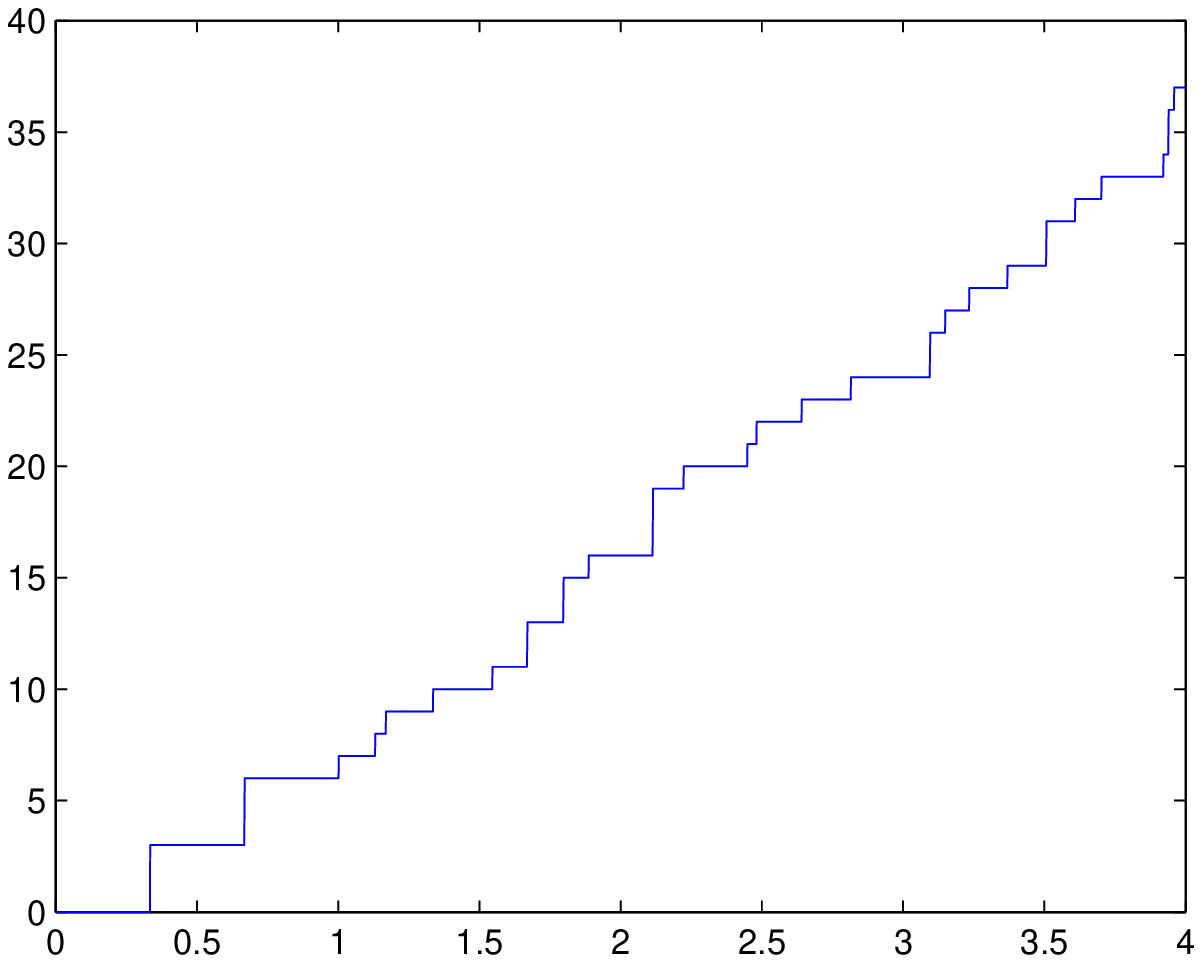}}
    \put(-65,0){\tiny Time}
    \put(-185,0){\tiny Time}
    \put(-240,58){\tiny $V$}
    \put(-124,62){\tiny $N_E$}
    \put(-104,5.5){\tiny $0$}
    \caption{Plots of (a) the evolution of the Lyapunov function $V$
      and (b) the total number $N_E$ of events of the
      \algoeventconsensus with $\sigma_i = 0.999$ for all~$i$. 
      The network consists of $5$ agents with
      communication topology described by the weight-balanced digraph
      $(\until{5}, \{(1,2), (2,3), (2,4), (3,4), (4,5), (5,1),
      (5,2)\})$ with weights $(1,1,0.5,1,1.5,1,0.5)$.  The initial
      condition is $x(0) = [-1, 0, 2, 2, 1]^T$.  }\label{fig:sim2}
  \end{figure}
}

We have also compared the \algoperiodic with a periodic implementation
of Laplacian consensus, cf.~\citep{ROS-JAF-RMM:07}. For the latter,
trajectories are guaranteed to converge if the timestep satisfies~$h <
{1}/{d^{\max}}$, where $d^{\max}$ is the maximum degree of the
graph~$\GG$. Figure~\ref{fig:sim3} shows this comparison using~$h =
0.1$ and also demonstrates the effect of $\{\sigma_i\}_{i=1}^N$ on the
executions of the \algoperiodic. For simplicity, we have used
$\sigma_i = \sigma$ to be the same for all agents in each execution.
One can observe the trade-off between communication and convergence
rate for varying~$\sigma$: higher~$\sigma$ results in less
communication but slower convergence compared to smaller values
of~$\sigma$.

{ \psfrag{oldwayonetwothree}[ll][ll]{\tiny \citep{EG-YC-HY-PA-DC:13}}%
  \psfrag{newwayonetwothree}[ll][ll]{\tiny Proposed}
  \psfrag{0}[cc][cc]{}
  \psfrag{1}[cc][cc]{\tiny $1$}%
  \psfrag{2}[cc][cc]{\tiny $2$} \psfrag{3}[cc][cc]{\tiny $3$}%
  \psfrag{4}[cc][cc]{\tiny $4$} \psfrag{5}[cc][cc]{\tiny $5$}%
  \psfrag{6}[cc][cc]{\tiny $6$} \psfrag{7}[cc][cc]{\tiny $7$}%
  \psfrag{8}[cc][cc]{\tiny $8$} \psfrag{10}[cc][cc]{\hspace*{-2mm}\tiny $10$}
  \psfrag{12}[cc][cc]{\tiny $12$} \psfrag{14}[cc][cc]{\tiny $14$}
  \psfrag{16}[cc][cc]{\tiny $16$}
  \psfrag{3.5}[cc][cc]{}
  \psfrag{2.5}[cc][cc]{}
  \psfrag{1.5}[cc][cc]{}
  \psfrag{0.5}[cc][cc]{}
  \psfrag{5}[cc][cc]{\hspace*{-2mm}\tiny $5$}
  \psfrag{20}[cc][cc]{\hspace*{-2mm}\tiny $20$}
  \psfrag{40}[cc][cc]{\hspace*{-2mm}\tiny $40$}
  \psfrag{15}[cc][cc]{\hspace*{-2mm}\tiny $15$}
  \psfrag{35}[cc][cc]{\hspace*{-2mm}\tiny $35$}
  \psfrag{30}[cc][cc]{\hspace*{-2mm}\tiny $30$}
  \psfrag{25}[cc][cc]{\hspace*{-2mm}\tiny $25$}
  \psfrag{150}[cc][cc]{\hspace*{-2mm}\tiny $150$}
  \psfrag{50}[cc][cc]{\hspace*{-2mm}\tiny $50$}
  \psfrag{100}[cc][cc]{\hspace*{-2mm}\tiny $100$}
  \begin{figure}[htb]
    \centering
    \subfigure[]{\includegraphics[width=.49\linewidth]{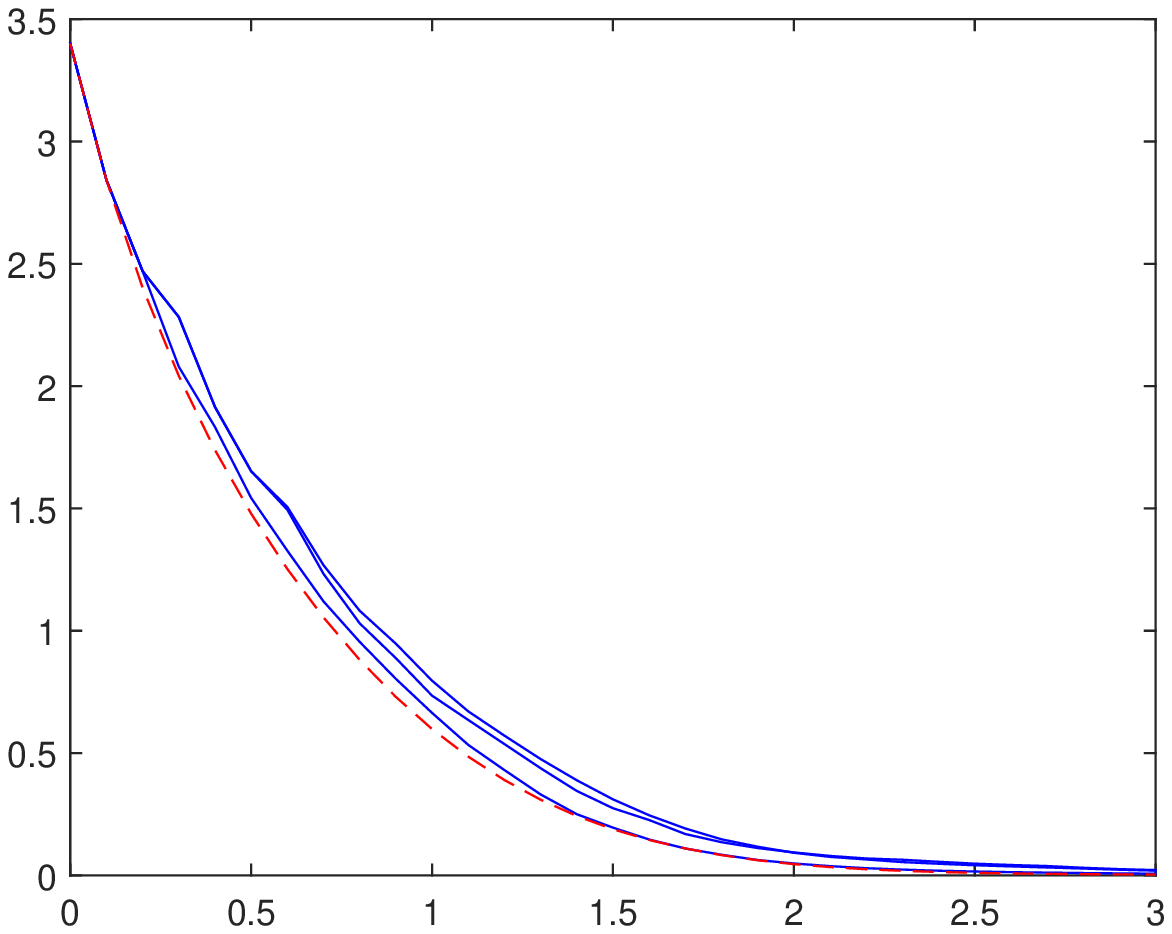}}
    \put(-104,5.5){\tiny $0$}
    \subfigure[]{\includegraphics[width=.49\linewidth]{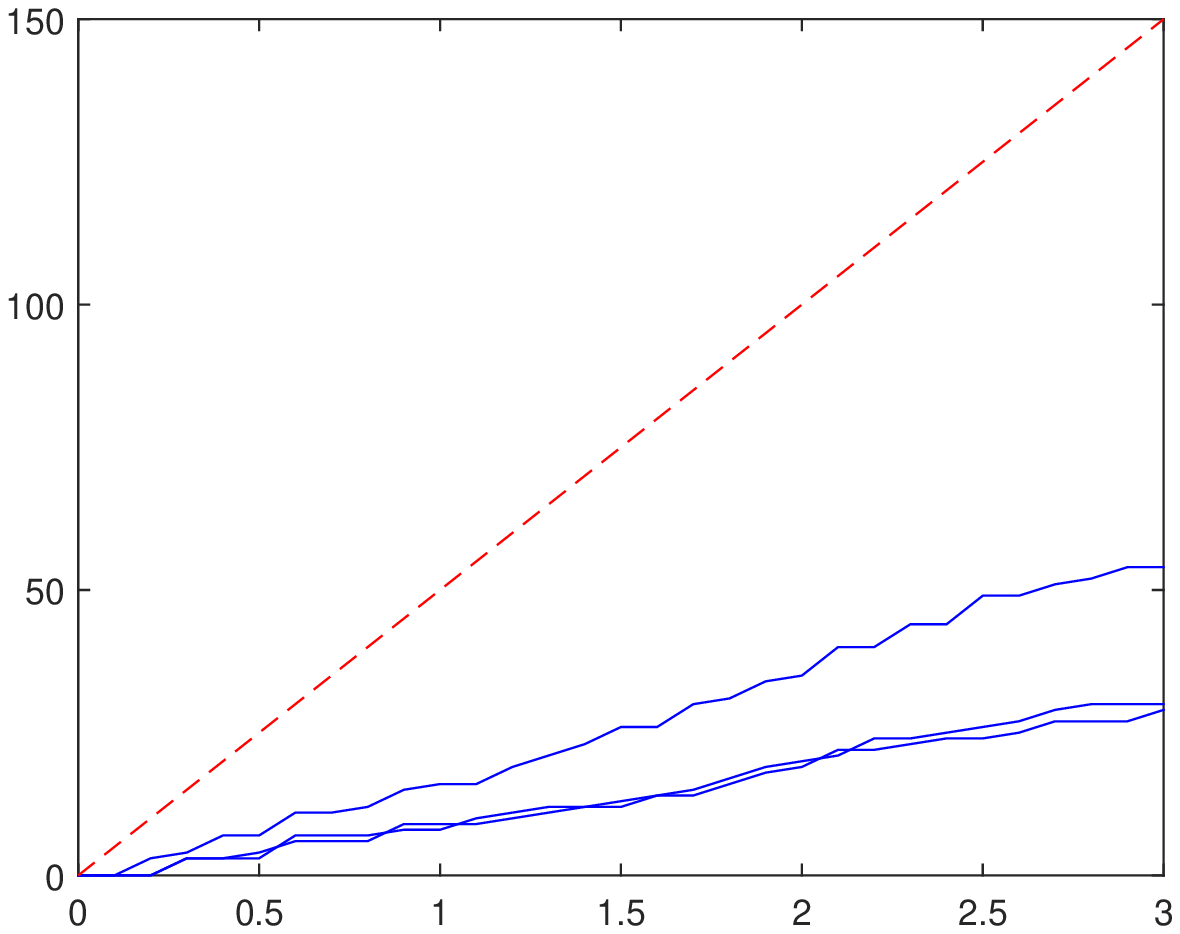}}
    \put(-65,0){\tiny Time}
    \put(-185,0){\tiny Time}
    \put(-235,60){\tiny $V$}
    \put(-121,65){\tiny $N_E$}
    \put(-104,5.5){\tiny $0$}
    \put(-45,35){\tiny $\sigma = 0.2$}
    \put(-36,25){\tiny $\sigma = 0.5$}
    \put(-36,16){\tiny $\sigma = 0.8$}
    \caption{Plots of (a) the evolution of the Lyapunov function $V$
      and (b) the total number $N_E$ of events of the \algoperiodic
      (with varying $\sigma$, solid blue) and a standard periodic
      Laplacian consensus algorithm (dashed red) with timestep $h =
      0.1$.  Network and initial condition are as in
      Figure~\ref{fig:sim2}.  }\label{fig:sim3}
  \end{figure}
}

\section{Conclusions}\label{se:conclusions}

We have proposed novel event-triggered communication and control
strategies for the multi-agent average consensus problem. Among the
novelties of our first design, we highlight that it works over
weight-balanced directed communication topologies, does not require
individual agents to continuously access information about the states
of their neighbors, and does not necessitate a priori agent knowledge
of global network parameters to execute the algorithm. We have shown
that our algorithms exclude the possibility of Zeno behavior and
identified conditions such that the network state exponentially
converges to agreement on the initial average of the agents' state.
We have also provided a lower bound on the convergence rate and
characterized the network convergence when the topology is switching
under a weaker form of connectivity.  Finally, we have developed a
periodic implementation of our event-triggered law that relaxes the
need for agents to evaluate the relevant triggering functions
continuously and provided a sufficient condition on the sampling
period that guarantee its the asymptotic correctness.  Future work
will explore scenarios with more general dynamics and physical sources
of error such as communication delays or packet drops, the extension
of our design and results to distributed convex optimization and other
coordination tasks, and further analysis of trigger designs that rule
out the possibility of Zeno behavior.

\section*{Acknowledgments}
This research was supported in part by NSF award CNS-1329619.

\end{document}